\documentclass{article}
\usepackage{latexsym,a4,amsmath,amsfonts}
\makeatletter
\@addtoreset{equation}{subsection}

\makeatother

\newcommand{\separate}{\medskip\noindent}

\def\newsection{ \separate
   \refstepcounter{subsection} 
   {\large\bf \thesubsection\kern.3em}
}
\def\mytheorem#1{
   \separate{\large\bf Theorem#1:\kern.3em}} 

\def\mylemma#1{
   \separate{\large\bf Lemma#1:\kern.3em}} 

\def\mycorollary#1{
   \separate{\large\bf Corollary#1:\kern.3em}} 

\def\myprop#1{
   \separate{\large\bf Proposition#1:\kern.3em}} 

\def\myremark#1{
   \separate{\large\bf Remark#1:\kern.3em}} 

\def\mydefinition#1{
   \separate{\large\bf Definition#1:\kern.3em}} 

\def\myconjecture#1{
   \separate{\large\bf Conjecture#1:\kern.3em}} 

\def\Proof{\separate\underline{Proof:}\kern1em}

\newcommand{\field}[1]{\mathbb{#1}}
\newcommand{\C}{\field{C}}
\newcommand{\R}{\field{R}}

\newcommand{\N}{\field{N}}
\newcommand{\Z}{\field{Z}}
\newcommand{\CPE}{{\mathchoice{\C{\rm P}_1}{\C{\rm P}_1}{\C\!\!
\mbox{\rm\tiny P}_1}{\C\!\!\mbox{\rm\tiny P}_1}}}

\def\QED{\hfill$\Box$}
\def\inv{^{-1}}
\def\bref#1{(\ref{#1})}
\def\tmatrix#1#2#3#4{
    \left(\begin{array}{cc} #1 & #2 \\ #3 & #4 \end{array}\right)}
\def\LieSL{{\mbox{\bf SL}}}

\def\LieSU{{\mbox{\bf SU}}}
\def\LieB{{\mbox{\bf B}}}

\def\LieU{{\mbox{\bf U}}}
\def\Liesl{{\mbox{\bf sl}}}

\def\dprime{{\prime\prime}}
\def\tprime{{\prime\prime\prime}}

\def\vecsigma{{\sigma\kern-2.25mm\sigma}}

\def\FD{{\cal D}}
\def\FF{{\cal F}}

\def\FRL{{{\cal R}\kern1mm_\lambda}}
\def\FS{{\cal S}}

\def\diff{{\rm d}}

\def\diag{{\mbox{\rm diag}}}

\def\Aut{{\mbox{\rm Aut}}}
\def\Sym{{\mbox{\rm Sym}}}

\def\const{{\mbox{\rm const}}}

\def\hg{\hat{g}}
\def\hf{\hat{f}}
\def\hg{\hat{g}}
\def\hh{\hat{h}}

\def\hF{\hat{F}}

\def\aquer{{\overline{a}}}
\def\bquer{{\overline{b}}}

\begin{document}

\renewcommand{\thefootnote}{\fnsymbol{footnote}}

\begin{center}
{\LARGE Investigation and application of the dressing action 
on surfaces of constant mean curvature}

\vskip1cm
\begin{minipage}{6cm}
\begin{center}
J.~Dorfmeister\footnotemark[1]\\
Department of Mathematics\\ University of Kansas\\ Lawrence, KS 66045
\end{center}
\end{minipage}
\begin{minipage}{6cm}
\begin{center}
G.~Haak\footnotemark[2]\\ Fachbereich Mathematik\\ TU Berlin\\ D-10623 Berlin
\end{center}
\end{minipage}
\vspace{0.5cm}

\end{center}
\footnotetext[1]{partially supported by NSF Grant DMS-9705479}
\footnotetext[2]{supported by Sonderforschungsbereich 288}

\section{Introduction}\label{introduction}\message{[introduction]}
The dressing method for the generating of new solutions of integrable
systems was introduced in 1979 by Zakharov
and Shabat \cite{ZaSh:1}. It was soon traced back to a natural loop
group action on the solution space of integrable
systems~\cite{AdvMo:1,ReSe:1}. 
With the introduction of the theory of integrable systems into
geometry, most notably the theory of surfaces of
constant mean curvature (CMC surfaces)~\cite{PiSt:1,Bo:1}, 
the dressing method also entered the realm of differential geometry.

Let it serve as an illustration of the power of the dressing action
that, due to~\cite{DoWu:1} (see also~\cite{BuPe:1}), all CMC
immersions of finite type, in particular all CMC tori, are included in
the dressing orbit of the standard cylinder. This result was used
in~\cite{DoHa:3} to reproduce the classification of CMC tori given
in~\cite{PiSt:1} in terms of the dressing action.

Moreover, while the integrable systems methods in general only apply
to a certain class of CMC surfaces, those without umbilics, the
dressing action can easily be applied also to CMC surfaces with
umbilics. The latter was made possible by a general loop group
theoretic approach to CMC surfaces, the so called DPW
method~\cite{DoPeWu:1}.

Similar to the Weierstra\ss\ representation of minimal surfaces, the
DPW method starts with a holomorphic function $E$ and a meromorphic
function $f$, both defined on some open, simply connected subset $\FD$
of the complex plane, and constructs an $S^1$-family of isometric
conformal CMC immersions $\Psi_{\lambda\in S^1}:\FD\rightarrow\R^3$,
the associated family, from these data. As a first application of the
dressing action, in~\cite{DoHa:1} and \cite{Wu:5} it was shown that
the dressing action can be used to describe the set of admissible
input data $(E,f)$ for the DPW method.

While $E\diff z^2$ is simply the Hopf differential of the resulting
CMC immersion, the function $f$ has no such simple geometric
interpretation. The problem shows especially if one wants to construct
CMC immersions $\Psi$ which are invariant under a symmetry group
$\Gamma\subset\Aut\FD$ of biholomorphic automorphisms of $\FD$, 
i.e.\ $\Psi\circ\gamma=T\Psi$, where
$\gamma\in\Gamma$ and $T$ is a (proper) Euclidean motion in space.  It
is clear, that $E\diff z^2$ has to be automorphic w.r.t.\ $\Gamma$,
i.e. $\gamma^\ast(E\diff z^2)=E\diff z^2$ or
$E\circ\gamma=(\gamma')^{-2}E$. But the meromorphic function $f$
satisfies no such automorphicity conditions. Instead, as was shown by
the authors in~\cite{DoHa:2}, $f$ transforms by complicated dressing
transformations under $\Gamma$.  For further reference
see~\cite{DoHa:2}. Thus, in the DPW formulation the understanding of
compact and symmetric CMC immersions is intimately related to the
understanding of the dressing action.  This served as a strong
incentive to further investigate the dressing action, in particular on
the meromorphic data of the DPW method.

Clearly, the Hopf differential is invariant under dressing.
Wu~\cite{Wu:4} constructed a large set of algebraic invariants under
dressing and was able to find normalized representatives in
each dressing orbit~\cite{Wu:3}. These results are particularly
useful, since the dressing orbits, as the finite type results
indicate, are very large.

In this paper we want to describe one of the major properties of
dressing as a group action: the isotropy group of a CMC immersion
under the dressing action.

As the main result of this paper we will show in
Theorem~\ref{mainresult} that for immersions whose Hopf differential
has zeroes, the isotropy group is trivial. In other words, the
dressing action is simple on surfaces with umbilics.  This is
essentially different from the situation for CMC surfaces of finite
type. Since these are all contained in the dressing orbit of the
standard cylinder, it is easy to see, that for these surfaces there is
an infinite dimensional isotropy group. In fact, as was shown
in~\cite{DoHa:3}, the existence of this isotropy group can be viewed
as the major prerequisite for Krichever's method. 

Theorem~\ref{mainresult} leads to several applications for CMC
immersions with umbilics. In this paper we present two
'no-go'-theorems for CMC surfaces with umbilics.

As a first candidate for a non-simply connected CMC surface with
umbilics we investigate the case $\FD=D_1$, the open unit disk, and
$f\equiv C$, where $C$ is a complex constant. In this case $f$ is
clearly invariant under any group $\Gamma\subset\Aut(D_1)$.  We will
show in Section~\ref{D1automorphisms}, that as long as the Hopf
differential has umbilics, such an $f$ cannot produce a CMC immersion
which is invariant under a group $\Gamma$ of biholomorphic
automorphisms of $D_1$.

As a second candidate we look for symmetric surfaces with umbilics in
the dressing orbits of arbitrary automorphic DPW data, a very large
class of surfaces. In Section~\ref{nosyminautom}, it will be shown
that for such CMC immersions all members of the associated family
$\Psi_\lambda$, $\lambda\in S^1$, constructed from the same DPW data,
share the same symmetry in $\R^3$. In particular, if
$\Psi=\Psi_{\lambda=1}$ is e.g.~compact in $\R^3$ of genus $g\geq2$
then all surfaces $\Psi_\lambda$ are compact with the same Fuchsian
group $\Gamma$. 

In view of the case of CMC tori, in which at most a countable number
of associated surfaces is compact, this may indicate that there are no
symmetric examples in the dressing orbits of automorphic DPW
data. However, as we want to emphasize, the latter is not a theorem,
just a conjecture.  In any case, the results of
Section~\ref{nosyminautom} show that even if there exists an
algebro-geometric method for CMC surfaces with umbilics, it will have
to look essentialy different from Krichever's method for the finite
type case.

At the end of this introduction we also would like to add, that in
spite of the discouraging results above, we were able to construct a
large family of non-simply connected CMC surfaces with umbilics using
the dressing action in the DPW method. These results will be published
elsewhere.

The authors want to thank F.~Pedit for many fruitfull discussions on
the subject which finally lead to Theorem~\ref{mainresult}.

\section{The dressing action}\refstepcounter{subsection} 
\label{dressingaction}
\message{[dressingaction]}
In this chapter we will review the definition and basic properties of
the dressing action on CMC surfaces in the framework of the DPW
construction. We will also investigate the isotropy group of a CMC
surface, represented by its meromorphic DPW data, under the dressing action.
For further reference see~\cite[Section~2]{DoHa:3}.

\newsection \label{loopgroups}
For each real constant $r$, $0<r<1$, 
let $\Lambda_r\LieSL(2,\C)_\sigma$ denote the group of 
smooth maps $g(\lambda)$ from $C_r$, the circle of radius $r$, to
$\LieSL(2,\C)$, which satisfy the twisting condition
\begin{equation} \label{DPWtwistcond}
g(-\lambda)=\sigma(g(\lambda)),
\end{equation}
where $\sigma:\LieSL(2,\C)\rightarrow\LieSL(2,\C)$ is defined by
conjugation with the Pauli matrix $\sigma_3=\tmatrix100{-1}$.
The Lie algebras of these groups, which we denote by
$\Lambda_r\Liesl(2,\C)_\sigma$, consist of maps 
$x:C_r\rightarrow\Liesl(2,\C)$, which satisfy a similar
twisting condition as the group elements
\begin{equation}
x(-\lambda)=\sigma_3 x(\lambda)\sigma_3.
\end{equation}
In order to make these loop groups complex Banach Lie groups, we equip them,
as in \cite{DoPeWu:1}, with some $H^s$-topology for $s>\frac12$.

Furthermore, we will use the following subgroups of
$\Lambda_r\LieSL(2,\C)_\sigma$: 
Let $\LieB$ be a subgroup of $\LieSL(2,\C)$ and 
$\Lambda_{r,B}^+\LieSL(2,\C)_\sigma$ be the group of maps in
$\Lambda_r\LieSL(2,\C)_\sigma$, which can be extended to holomorphic maps on
\begin{equation}
I^{(r)}=\{\lambda\in\C; |\lambda|<r\},
\end{equation}
the interior of the circle $C_r$, and take values in $\LieB$ at $\lambda=0$.
Analogously, let $\Lambda_{r,B}^-\LieSL(2,\C)_\sigma$ 
be the group of maps in $\Lambda_r\LieSL(2,\C)_\sigma$, which can be extended
to the exterior 
\begin{equation}
E^{(r)}=\{\lambda\in\CPE;|\lambda|>r\}
\end{equation}
of $C_r$ and take values in $\LieB$ at $\lambda=\infty$. 
If $\LieB=\{I\}$ (based loops) 
we write the subscript $\ast$ instead of $\LieB$, if
$\LieB=\LieSL(2,\C)$ we omit the subscript for $\Lambda$ entirely.

Also, by an abuse of notation, we will denote by 
$\Lambda_r\LieSU(2)_\sigma$ the subgroup of maps in
$\Lambda_r\LieSL(2,\C)_\sigma$, which can be extended holomorphically
to the open annulus
\begin{equation}
A^{(r)}=\{\lambda\in\C;r<|\lambda|<\frac{1}{r}\}
\end{equation}
and take values in $\LieSU(2)$ on the unit circle.

Corresponding to these subgroups, we analogously define Lie subalgebras of
$\Lambda_r\Liesl(2,\C)_\sigma$.

We quote the following results from~\cite{Mc:1} and \cite{DoPeWu:1}:

\separate (i) For each solvable subgroup $\LieB$ of 
$\LieSL(2,\C)$, which satisfies $\LieSU(2)\cdot\LieB=\LieSL(2,\C)$ and
$\LieSU(2)\cap\LieB=\{I\}$, multiplication 
$$
\Lambda_r\LieSU(2)_\sigma\times\Lambda_{r,B}^+\LieSL(2,\C)_\sigma
\longrightarrow\Lambda_r\LieSL(2,\C)_\sigma
$$
is a diffeomorphism onto.
The associated splitting
\begin{equation} \label{Iwasawa}
g=F g_+
\end{equation}
of an element $g$ of $\Lambda_r\LieSL(2,\C)_\sigma$, s.t.\
$F\in\Lambda_r\LieSU(2)_\sigma$ and
$g_+\in\Lambda^+_{r,B}\LieSL(2,\C)_\sigma$
will be called Iwasawa decomposition.

\separate (ii) Multiplication 
\begin{equation} \label{Birkhoff}
\Lambda^-_{r,\ast}\LieSL(2,\C)_\sigma\times\Lambda_r^+\LieSL(2,\C)_\sigma
\longrightarrow\Lambda_r\LieSL(2,\C)_\sigma
\end{equation}
is a diffeomorphism onto the open and dense subset 
$\Lambda^-_{r,\ast}\LieSL(2,\C)_\sigma\cdot\Lambda_r^+\LieSL(2,\C)_\sigma$
of $\Lambda_r\LieSL(2,\C)_\sigma$, called the ``big cell'' \cite{SeWi:1}.
The associated splitting
\begin{equation}
g=g_-g_+
\end{equation}
of an element $g$ of the big cell, where
$g_-\in\Lambda^-_{r,\ast}\LieSL(2,\C)_\sigma$ and
$g_+\in\Lambda_r^+\LieSL(2,\C)_\sigma$, 
will be called Birkhoff factorization.

\newsection \label{dressingdef}
Let $\Psi:\FD\rightarrow\R^3$ be a conformal CMC-immersion. Define the
extended frame $F(z,\lambda):\FD\rightarrow\Lambda\LieSU(2)_\sigma$ as
in~\cite{DoPeWu:1} (see also the appendix of~\cite{DoHa:1}).
Furthermore, define $g_-:\FD\rightarrow\Lambda^-_{r,\ast}\LieSL(2,\C)_\sigma$
by the Birkhoff splitting
\begin{equation} \label{gmindef}
F(z,\lambda)=g_-(z,\lambda)g_+(z,\lambda).
\end{equation}
Then $g_-$ is a meromorphic function on $\FD$ with poles in the set
$\FS\subset\FD$ of points, where $F(z,\lambda)$ is not in the ``big
cell'', i.e.\ where the Birkhoff splitting~\bref{gmindef} of 
$F(z,\lambda)$ is not defined. It should also be noted that,
by~\cite[Lemma~2.2]{DoHa:3}, the maximal analytic continuation of
$g_-$ does not depend on the chosen radius $r$. I.e.\ the
meromorphic potential of a CMC immersion does not depend on $r$.

For given meromorphic $g_-$, we can recover the extended frame
$F$ by the Iwasawa decomposition
\begin{equation}
g_-=Fg_+\inv.
\end{equation}
For smoothness questions, see~\cite{DoHa:1}.

Next, we define the dressing action of
$\Lambda_r^+\LieSL(2,\C)_\sigma$, $0<r\leq1$, on
$\FF$, the set of extended frames of
CMC-immersions. For $B(z,\lambda)\in\FF$ and
$h_+\in\Lambda_r^+\LieSL(2,\C)_\sigma$ we set
\begin{equation} \label{dressing}
h_+(\lambda)F(z,\lambda)=(h_+.F)(z,\lambda)q_+(z,\lambda),
\end{equation}
where the r.h.s.\ of~\bref{dressing} is defined by the Iwasawa
decomposition in $\Lambda_r\LieSL(2,\C)_\sigma$ of $h_+F$, i.e.\
$q_+:\FD\rightarrow\Lambda^+_r\LieSL(2,\C)_\sigma$. In
addition at $\lambda=0$
the matrix $q_+(z,\lambda)$ takes values in the solvable subgroup $B$ of
$\LieSL(2,\C)$, s.t.\ 
\begin{equation} \label{trivialintersection}
B\cap\LieSU(2)=\{I\}.
\end{equation}
It is easily proved (see e.g.~\cite{BuPe:1})
that $h_+.F$ is again in $\FF$. Therefore,
Eq.~\bref{dressing} really defines an action on $\FF$.
On the matrices $g_-$ defined by~\bref{gmindef} the dressing is
defined by
\begin{equation} \label{gmindressing}
h_+(\lambda)g_-(z,\lambda)=\hg_-(z,\lambda)p_+(z,\lambda).
\end{equation}
Here, $\hg_-=h_+.g_-$ and
$p_+:\FD\rightarrow\Lambda_r^+\LieSL(2,\C)_\sigma$ are defined by the Birkhoff
splitting~\bref{gmindef} of $h_+g_-$. Note that $h_+.F=\hg_-\hg_+$ for
some $\hg:\FD\rightarrow\Lambda_r^+\LieSL(2,\C)_\sigma$. Since $g_-$ and
$\hg_-$ are both meromorphic in $z$, also
$p_+=\hg_+q_+g_+\inv$ is meromorphic in $z$.

The extended frames are normalized by
\begin{equation} \label{Finitial}
F(0,\lambda)=I,\kern1cm\lambda\in S^1,
\end{equation}
which implies
\begin{equation} \label{gmininitial}
g_-(0,\lambda)=I,\kern1cm\lambda\in S^1.
\end{equation}
Let now the meromorphic potential be defined by 
\begin{equation}
\xi(z,\lambda)=g_-\inv\diff g_-,
\end{equation}
then it is of the form
\begin{equation} \label{xiform}
\xi(z,\lambda)=\lambda\inv\tmatrix0f{\frac{E}{f}}0\diff z,
\end{equation}
where $f$ is a nonvanishing meromorphic function.
We will always assume $E\not\equiv0$, 
i.e.\ we will exclude the case
that the surface is part of a round sphere.
To construct a CMC-immersion from a given meromorphic potential of the
form~\bref{xiform}, the functions $f$ and $E$ cannot be
chosen arbitrarily. They have to satisfy additional conditions, given in
\cite{DoHa:1}.

The matrix $g_-$ and therefore also the frame $F$ are uniquely 
determined by the meromorphic potential and the initial
condition~\bref{Finitial}.

From~Eq.~\bref{gmindressing} it follows, that $\xi$ transforms under
dressing with $h_+\in\Lambda_r^+\LieSL(2,\C)_\sigma$ as
\begin{equation} \label{xidressing}
h_+.\xi=p_+\inv\xi p_++p_+\inv\diff p_+=\lambda\inv
\tmatrix0{h_+.f}{\frac{E}{h_+.f}}0\diff z.
\end{equation}
Note, that $E\diff z^2$, the Hopf differential of the CMC-immersion
$\Psi$, is invariant under dressing.

We now set
\begin{equation} \label{abcddef}
p_+=\tmatrix abcd.
\end{equation}
Then, for the matrix entries of $p_+$ we get with $\hf=h_+.f$:
\begin{eqnarray}
\lambda a^\prime & = & b\frac{E}{\hf}-fc,\label{genaprime}\\
\lambda b^\prime & = & a\hf-fd,\label{genbprime}\\
\lambda c^\prime & = & d\frac{E}{\hf}-\frac{E}{f}a,\label{gencprime}\\
\lambda d^\prime & = & c\hf-\frac{E}{f}b,\label{gendprime}
\end{eqnarray}
where $(\cdot)'$ denotes differentiation w.r.t.\ $z$.

\newsection \label{isotropy}
Let us investigate the isotropy group $I(F)$ of an extended frame $F$
under dressing. For an extended frame $F$,
$I(F)$ is defined as the group of all
$h_+\in\Lambda_r^+\LieSL(2,\C)_\sigma$, s.t.\ $h_+.F=F$.

Assume, that $h_+\in I(F)$. Then Eq.~\bref{gmindressing} defines a
meromorphic function $p_+$ on $\FD$, s.t.\
\begin{equation} \label{xiiso}
\xi=h_+.\xi=p_+\inv\xi p_++p_+\inv\diff p_+,
\end{equation}
or, using~\bref{genaprime}--\bref{gendprime},
\begin{eqnarray}
\lambda a^\prime & = & b\frac{E}{f}-fc=-\lambda d^\prime,\label{adprime}\\
\lambda b^\prime & = & (a-d)f,\label{bprime}\\
\lambda c^\prime & = & (d-a)\frac{E}{f}.\label{cprime}
\end{eqnarray}

\mylemma{} {\em
Let $F$ be an extended frame and let $\xi$ be the associated
meromorphic potential.
Let $h_+\in I(F)$ and let
$p_+:\FD\rightarrow\Lambda_r^+\LieSL(2,\C)_\sigma$
be the associated solution of
Eq.~\bref{xiiso}. Define $b$ as the upper right entry of $p_+$.
If $b\equiv0$, then $h_+=I$.
}

\separate\underline{Proof:}
If $b\equiv0$, then Eq.~\bref{bprime} gives $a=d$.
This together with $\det p_+=ad=1$ shows, that
$a=d\equiv\pm1$. By~\bref{adprime} we get $c\equiv0$.
This implies, that $p_+=\pm I$. Therefore, by Eq.~\bref{gmindressing},
\begin{equation}
h_+g_-=\pm g_-.
\end{equation}
This shows, that $h_+=\pm I$. The case $h_+=-I$ is excluded, since in
this case, by Eq.~\bref{dressing}, $q_+=-I\in B\cap\LieSU(2)$, which
contradicts~\bref{trivialintersection}.
\QED

\separate
We want to investigate the set of scalar meromorphic functions
$(a,b,c,d)$ on $\FD$,
s.t.\ the set of equations~\bref{adprime}--\bref{cprime} is
satisfied. To this end, we first rewrite~\bref{adprime}--\bref{cprime}
as a single third order differential equation in $b$:

First, we differentiate Eq.~\bref{bprime} and use Eq.~\bref{adprime} to
get
\begin{equation} \label{bdprime}
\lambda b^\dprime=\frac{f^\prime}{f}\lambda b^\prime+2fa^\prime
\end{equation}
or
\begin{equation} \label{aprime2}
a^\prime=\frac{\lambda}{2f}\left(b^\dprime-\frac{f^\prime}{f}b^\prime\right).
\end{equation}
From Eq.~\bref{adprime} we also get, using~\bref{aprime2},
\begin{equation} \label{ceq}
c=-\frac{\lambda}{f}a^\prime+\frac{E}{f^2}b=-\frac{\lambda^2}{2f^2}b^\dprime
+\frac{\lambda^2f^\prime}{2f^3}b^\prime+\frac{E}{f^2}b.
\end{equation}
From~\bref{bprime} and~\bref{cprime} it follows, that
\begin{equation} \label{bcprime}
b^\prime\frac{E}{f}=-c^\prime f.
\end{equation}
By differentiating Eq.~\bref{adprime} we furthermore get, using
Eq.~\bref{bcprime},
\begin{equation} \label{adprime2}
\lambda a^\dprime=2b^\prime\frac{E}{f}-f^\prime c
+b\left(\frac{E}{f}\right)^\prime.
\end{equation}
Differentiating~\bref{bdprime} and using~\bref{aprime2}, \bref{adprime2}
and~\bref{ceq}, we get the following ordinary differential equation in
$b$:
\begin{equation} \label{bODE}
\lambda^2\left(b^\tprime-3\frac{f^\prime}{f}b^\dprime
-\left(\left(\frac{f^\prime}{f}\right)^\prime
-2\left(\frac{f^\prime}{f}\right)^2\right)b^\prime\right)
=4Eb^\prime+2\left(E^\prime-2\frac{f^\prime}{f}E\right)b.
\end{equation}
Since $p_+$ takes values in the twisted loop group
$\Lambda_r\LieSL(2,\C)_\sigma$, the function $b(z,\lambda)$ is odd in
$\lambda$. If we write
\begin{equation} \label{bdecomp}
b(z,\lambda)=\sum_{n=0}^\infty b_n(z)\lambda^n,
\end{equation}
then all coefficients $b_n$ for which $n$ is even are identically
zero.

From Eq.~\bref{bODE} we get
\begin{equation} \label{b1eq}
2Eb_1^\prime+\left(E^\prime-2\frac{f^\prime}{f}E\right)b_1=0
\end{equation}
and the recursion relation
\begin{equation} \label{bneq}
b_{n-2}^\tprime-3\frac{f^\prime}{f}b_{n-2}^\dprime
-\left(\left(\frac{f^\prime}{f}\right)^\prime
-2\left(\frac{f^\prime}{f}\right)^2\right)b_{n-2}^\prime
=4Eb_n^\prime+2\left(E^\prime-2\frac{f^\prime}{f}E\right)b_n,\;n\geq3.
\end{equation}
We can solve Eq.~\bref{b1eq} for $b_1$ and we get
\begin{equation} \label{b1solution}
b_1=C_b\sqrt{\frac{f^2}{E}},
\end{equation}
where $C_b$ is a constant.

The conjugation of
Eq.~\bref{xidressing} with the matrix $\tmatrix0110$ leads to a
transformation of the system of
equations~\bref{adprime}--\bref{cprime}.
In terms of matrix entries this transformation reads
\begin{equation}
f\longrightarrow\frac{E}{f},\kern1cm a\longleftrightarrow d,\kern1cm
b\longleftrightarrow c.
\end{equation}
Using this transformation, we can immediately write down the ODE for
the coefficient $c$ by replacing $f$ by $\frac{E}{f}$ in
Eq.~\bref{bODE}:
\begin{equation} \label{cODE}
\lambda^2\left(c^\tprime
-3\left(\frac{E^\prime}{E}-\frac{f^\prime}{f}\right)c^\dprime
-\left(\left(\frac{E^\prime}{E}-\frac{f^\prime}{f}\right)^\prime
-2\left(\frac{E^\prime}{E}-\frac{f^\prime}{f}\right)^2\right)c^\prime\right)
=4Ec^\prime-2\left(E^\prime-2\frac{f^\prime}{f}E\right)c.
\end{equation}
We define the $\lambda$-coefficients $c_n$ of $c$ by
\begin{equation} \label{cdecomp}
c(z,\lambda)=\sum_{n=0}^\infty c_n(z)\lambda^n.
\end{equation}
Like $b$ also $c$ is odd in $\lambda$. Therefore, for $n$ even, the
coefficients $c_n$ vanish.
From~\bref{cODE} we get a recursion relation for the $c_n$:
\begin{equation} \label{c1eq}
2Ec_1^\prime-\left(E^\prime-2\frac{f^\prime}{f}\right)c_1=0
\end{equation}
and
\begin{equation} \label{cneq}
c_{n-2}^\tprime
-3\left(\frac{E^\prime}{E}-\frac{f^\prime}{f}\right)c_{n-2}^\dprime
-\left(\left(\frac{E^\prime}{E}-\frac{f^\prime}{f}\right)^\prime
-2\left(\frac{E^\prime}{E}-\frac{f^\prime}{f}\right)^2\right)c_{n-2}^\prime
=4Ec_n^\prime-2\left(E^\prime-2\frac{f^\prime}{f}E\right)c_n,
\end{equation}
$n\geq3$. We can solve Eq.~\bref{c1eq} for $c_1$ and we get
\begin{equation} \label{c1solution}
c_1=C_c\sqrt{\frac{E}{f^2}},
\end{equation}
where $C_c$ is a constant. 

Using only Eq.~\bref{bneq} we get

\mytheorem{} {\em
Let $\Psi:\FD\rightarrow\R^3$ be a conformal CMC-immersion with
extended frame $F(z,\lambda)$ and define the dressing action as above.
If the isotropy group $I(F)$ of $F$ under dressing is nontrivial, i.e.\
$I(F)\neq\{I\}$, then the Hopf differential $E$ of $\Psi$ is the square
of a meromorphic function.
}

\separate\underline{Proof:}
We assume that there exists $I\neq h_+\in I(F)$. 
Let $p_+$ be the associated solution
of~\bref{xiiso}. We define $a$, $b$, $c$, and
$d$ as in Eq.~\bref{abcddef}. Define $b_n$, $n\in\N$,
by~\bref{bdecomp}. By Lemma~\ref{isotropy}, $b\not\equiv0$.
Thus there exists a smallest index $N\in\N$, $N$ odd,
for which $b_n\not\equiv0$, i.e.\ $b_N\not\equiv0$ and $b_n\equiv0$
for $n<N$. The $\lambda^N$-coefficient
$b_N$ is meromorphic and satisfies
\begin{equation} \label{bNeq}
2Eb_N^\prime+\left(E^\prime-2\frac{f^\prime}{f}E\right)b_N=0,
\end{equation}
since the l.h.s.\ of Eq.~\bref{bneq} vanishes.
Eq.~\bref{bNeq} has the solution
\begin{equation}
b_N=C\sqrt{\frac{f^2}{E}},
\end{equation}
where $C$ is a complex constant, $C\neq0$. Since $b_N$ is meromorphic,
we get that the Hopf differential
\begin{equation}
E=\left(C\frac{f}{b_N}\right)^2
\end{equation}
is the square of a meromorphic function.
\QED

\mycorollary{} {\em
Let $\Psi$ and $F$ be as in Theorem~\ref{isotropy}.
If the Hopf differential of
$\Psi$ has a zero of odd order, then $I(F)=\{I\}$.
}

\newsection \label{mainresult}
We can actually extend Corollary~\ref{isotropy} to all surfaces with umbilics.
Let us state the main result of this paper:

\mytheorem{} {\em
Let $\Psi:\FD\rightarrow\R^3$ be a conformal CMC-immersion with
extended frame $F(z,\lambda)$ and define the dressing action as above.
If the surface defined by $\Psi$ has an umbilic, i.e.\ if its Hopf
differential $E$ has a zero,
then the isotropy group $I(F)$ of $F$ under dressing is trivial,
i.e.\ $I(F)=\{I\}$.
}

\separate
Before we prove Theorem~\ref{mainresult},
we will draw a simple conclusion from Eq.~\bref{bneq} and Eq.~\bref{cneq}:

\mylemma{} {\em
Let $\xi=\lambda\inv\tmatrix0f{\frac{E}{f}}0\diff z$ 
be a meromorphic potential and let $p_+$ be a meromorphic
matrix function which satisfies Eq.~\bref{xiiso}. Define
$b(z,\lambda)$ and $c(z,\lambda)$
by Eq.~\bref{abcddef}. Then the following holds:

1. If $f$ is defined at $z_0\in\FD$ then $b$ is defined at $z_0$. 

2. If $f$ has a pole of order $j$ at $z_0$, then $b$ has at most a pole of
order $2(j-1)$ at $z_0$.

3. If $\frac{E}{f}$ is defined at $z_0$, then $c$ is defined at $z_0$.

4. If $\frac{E}{f}$ has a pole of order $j$ at $z_0$, then $c$ has at
most a pole of order $2(j-1)$ at $z_0$.
}

\separate\underline{Proof:}
We know, that the coefficients $b_n(z)$ and $c_n(z)$ are all meromorphic in
$\FD$. 
Let us denote by $k_n$ the order of the pole of $b_n$ at $z_0$ for
$n\in\N$, $k_n=0$ if $b_n$ is defined at $z_0$.
The function $\frac{f^\prime}{f}$ has at most a simple pole at $z_0$.
Let $n_f\in\Z$ be the residue of $\frac{f^\prime}{f}$ at $z_0$.

If $k_{n-2}>0$ then the l.h.s.\ of Eq.~\bref{bneq} has at most a pole
of order $k_{n-2}+3$ at $z_0$. The function $b_{n-2}$ is of the form
\begin{equation}
b_{n-2}(z)=\beta z^{-k_{n-2}}+v(z)
\end{equation}
with $\beta\neq0$, and $z^{k_{n-2}-1}v(z)$ locally holomorphic at $z_0$.
The coefficient of $z^{-(k_{n-2}+3)}$ on the l.h.s.\
of Eq.~\bref{bneq} is given by
\begin{equation}
\beta(-(k_{n-2}+2)(k_{n-2}+1)k_{n-2}-3n_f(k_{n-2}+1)k_{n-2}
+(-n_f-2n_f^2)k_{n-2}).
\end{equation}
Therefore, the l.h.s.\ of Eq.~\bref{bneq} has a pole of order
$k_{n-2}+3>3$ at $z_0$ iff
\begin{equation} \label{degcond}
(k_{n-2}+2)(k_{n-2}+1)+3n_f(k_{n-2}+1)+n_f+2n_f^2\neq0.
\end{equation}
This in turn is equivalent to
\begin{equation} \label{degcond2}
k_{n-2}\neq-2(n_f+1)\kern2cm\mbox{\rm and}\kern2cm k_{n-2}\neq-(n_f+1).
\end{equation}

1. If $n_f\geq0$, i.e.\ $f$ has no pole at $z_0$,
then the condition~\bref{degcond2} is always satisfied if $k_{n-2}>0$.
Therefore, if $b_{n-2}$ has a pole of order $k_{n-2}>0$ at
$z_0$, then the l.h.s.\ of 
Eq.~\bref{bneq} has a pole of order $k_{n-2}+3>3$ at $z_0$.
Thus, also the r.h.s.\ of Eq.~\bref{bneq} has
a pole of order $\geq4$ at $z_0$. This is only possible, if $b_n$ has a pole at
$z_0$. If $k_n>0$ then the r.h.s.\ of Eq.~\bref{bneq} has a pole of
order at most $k_n-m+1$ at $z_0$. Here, $m\geq0$ 
is the zero order of $E$ at $z_0$,
$m=0$ if $E(z_0)\neq0$. By comparing the pole orders, we get
\begin{equation} \label{poleeq}
k_n\geq k_{n-2}+m+2>k_{n-2}.
\end{equation}
Let us assume, that there exists $N\in\N$, s.t.\ $k_N>0$. 
Using~\bref{poleeq} we get that $b_{N+2l}$, $l>0$, has a pole of order
\begin{equation}
k_{N+2l}\geq k_N+l(m+2),\;l>0,
\end{equation}
at $z_0$. 
It follows, that $b$ has an essential singularity at $z=z_0$ for all
$\lambda\in S^1$. This contradicts the meromorphicity of
$p_+$. Therefore, all $b_n$ are holomorphic in $\FD$.

2. If $n_f<0$ then $j=-n_f$, and the
condition~\bref{degcond} is certainly satisfied if
\begin{equation} \label{niceeq}
k_{n-2}>-2(n_f+1).
\end{equation}
In this case, we can argue as in the proof of 1.\
that Eq.~\bref{poleeq} is satisfied for $k_n$. Since therefore
$k_n>k_{n-2}$, we get that also for $k_n$ Eq.~\bref{niceeq} is
satisfied. Assume, that there exists $N\in\N$, s.t.\
$k_N>-2(n_f+1)$. Then, by the argument above, Eq.~\bref{niceeq} is
satisfied for all $n\geq N$ and we get, as in the first part:
\begin{equation}
k_{N+2l}\geq k_N+l(m+2),\;l>0.
\end{equation}
This shows that $b$ has an essential singularity at $z_0$,
contradicting the meromorphicity of $p_+$. Therefore, $b$ can have at most
a pole of order $-2(n_f+1)=2(j-1)$.

3.\ and 4.\ follow from the proof of 1.\ and 2.\
by replacing $f$ by $\frac{E}{f}$ and Eq.~\bref{bneq}
by Eq.~\bref{cneq}.
\QED

\newsection \label{singularities}
The fact that dressing is a group action allows us to prove the following

\mylemma{} {\em
Let $F$ be the extended frame of a CMC immersion and let
$\hh_+\in\Lambda_r^+\LieSL(2,\C)_\sigma$. Define $\hF=\hh_+.F$. Then
the isotropy groups $I(F)$ and $I(\hF)$ are isomorphic. The
isomorphism is given by conjugation with $\hh_+$ in 
$\Lambda_r^+\LieSL(2,\C)_\sigma$.
}

\separate\underline{Proof:}
Assume, that $h_+\in I(F)$. Then we have
\begin{equation}
h_+.F=F.
\end{equation}
Since dressing is a group action, we can rewrite this as
\begin{equation}
h_+.(\hh_+\inv.\hF)=\hh_+\inv.\hF.
\end{equation}
This shows that
\begin{equation}
(\hh_+h_+\hh_+\inv).\hF=\hF,
\end{equation}
i.e.\ $\hh_+h_+\hh_+\inv\in I(\hF)$. Conversely, for each $h_+\in I(\hF)$
we see in the same way, that $\hh_+\inv h_+\hh_+\in I(F)$.
\QED

\separate 
Using Lemma~\ref{isotropy}, Lemma~\ref{mainresult}, and
Lemma~\ref{singularities}, we can now give the

\separate\underline{Proof of Theorem~\ref{mainresult}:} 
Let $f$ be defined by Eq.~\bref{xiform}. We assume, that $E$ has a
zero of order $m>0$ at some $z_0\in\FD$. Let $h_+\in I(F)$ and let
$p_+$ be the corresponding solution of~\bref{xiiso}. Let $b$ be the
upper right entry of $p_+$ and let $b_n$, $n\in\N$, be defined by
Eq.~\bref{bdecomp}. 

Assume $b\not\equiv0$. If $N$ is the smallest index for which
$b_n\not\equiv0$ then by the same argument as in the proof of
Theorem~\ref{isotropy} we get
\begin{equation} \label{bNeq2}
b_N=C\sqrt{\frac{f^2}{E}}.
\end{equation}

Case I: $f$ is locally holomorphic without zero at $z_0$. Then $b$ is, by
Lemma~\ref{mainresult}, defined and locally holomorphic at $z_0$.
But $b_N$ has by~\bref{bNeq2} a pole of order $m$ at $z_0$, a contradiction.
Therefore, we get $b\equiv0$ and, by Lemma~\ref{isotropy},
$h_+=I$. This shows that in this case $I(F)=\{I\}$.

Case II: $f$ has a pole or zero at $z_0$. In
\cite[Section~3.12]{DoHa:1} it was shown, that in
the $r=1$-dressing orbit of each extended frame $F$ there is a frame $\hF$,
s.t.\ the function $\hf$, defined by the associated
meromorphic potential, has neither a pole nor a zero at $z_0$. Since
$\Lambda^+\LieSL(2,\C)_\sigma$ is a subgroup of
$\Lambda_r^+\LieSL(2,\C)_\sigma$ for each $0<r\leq1$, the same
statement holds for each $r$-dressing orbit of an extended frame. The
isotropy group of $\hF$ is therefore, by the proof of Case I, trivial,
i.e.\ $I(\hF)=\{I\}$.
By Lemma~\ref{singularities}, for two elements $F$ and
$\hF$ in the same dressing orbit the isotropy groups $I(F)$ and
$I(\hF)$ are isomorphic. This shows that $I(F)=I(\hF)=\{I\}$, which
finishes the proof.
\QED

\separate Let us emphasize the main result again by reformulating it
in the following form:
 
\mycorollary{} {\em
If the isotropy group of a CMC immersion under dressing is nontrivial,
then the surface has no umbilics.
}

\section{Applications}\label{applications}\message{[applications]}

In this section we will apply the results of the last section, i.e.\
Theorem~\ref{mainresult} to symmetric surfaces of constant mean
curvature.

Let now $\xi$ be a meromorphic potential given by \bref{xiform} with
$f\equiv C=\const$ and $E$ a holomorphic function {\em with zeroes} on the open
unit disk $\FD=D_1$. I.e.\ the
CMC immersion $\Psi:\FD\rightarrow\R^3$ associated to $\xi$
has umbilics.

We will show, that for such a surface the symmetry group $\Sym(\Psi)$
of biholomorphic automorphisms which leave the surface invariant up to
a proper Euclidean motion in space, will never contain a fixed point
free automorphism. In other words, it is not possible using a constant
function $f$ to construct a not simply connected CMC surface over $\FD=D_1$
with umbilics.

\newsection \label{D1automorphisms}

The following is well known:

\mylemma{} {\em
The group of biholomorphic automorphisms $\Aut(D_1)$ of the unit disk
is the following group of Moebius transformations:
$$
\Aut(D_1)=\{\gamma:z\longmapsto \frac{az+b}{\bquer
z+\aquer},\;a,b\in\C, |a|^2-|b|^2=1\}.
$$
An element $\gamma\in\Aut(D_1)$ has a fixed point inside $D_1$
iff it describes a rotation around the origin $z=0$, 
i.e.\ $|a|=1$ and $b=0$.
}

We will use this result together with Theorem~\ref{mainresult} to
prove

\mytheorem{} {\em
Let $\Psi:D_1\rightarrow\R^3$ be a CMC immersion with umbilics
whose meromorphic potential is of the form
\begin{equation} \label{xifconst}
\xi=\lambda\inv\tmatrix0f{\frac{E}{f}}0\diff z,\;f\equiv C\in\C.
\end{equation}
If the extended frame $F:D_1\rightarrow\Lambda\LieSU(2)_\sigma$ of
$\Psi$ satisfies
\begin{equation} \label{Ftrafo}
F(\gamma(z),\lambda)=\chi(\lambda)F(z,\lambda)k(z)
\end{equation}
for some $\gamma\in\Aut(D_1)$, $\chi\in\Lambda\LieSU(2)_\sigma$, and
$k:D_1\rightarrow\LieU(1)$ then
$\gamma$ is a rotation around the origin $z=0\in D_1$.
}

\mycorollary{} {\em
It is impossible to obtain non-simply connected CMC immersions
$\Psi:D_1\rightarrow\R^3$ with
umbilics, in particular compact CMC surfaces of genus $g\geq2$, 
from meromorphic potentials of the form~\bref{xifconst}.
}

\separate\underline{Proof of Theorem~\ref{D1automorphisms}:}\kern1em
We first write~\bref{Ftrafo} as
\begin{equation}
(g_-g_+)\circ\gamma=\chi_-\chi_+g_-g_+k
\end{equation}
where $F=g_-g_+$ and $\chi=\chi_-\chi_+$ denote the Birkhoff splitting
of $F$ and $\chi$, respectively. Of course, $g_-\inv\diff g_-=\xi$ is
the meromorphic potential.
Since $\chi(\lambda)=F(\gamma(0),\lambda)$, we get
$$
\chi_-=g_-(\gamma(0)),\kern1cm\chi_+k=g_+(\gamma(0)).
$$
If we set
\begin{equation}
\hg_-=\chi_-\inv g_-\circ\gamma)=\chi_+g_-g_+k(g_+\circ\gamma)\inv
\end{equation}
then
\begin{equation}
\hg_-\inv\diff\hg_-=(g_-\circ\gamma)\inv\diff(g_-\circ\gamma)=\xi\circ\gamma.
\end{equation}
Since also $\hg_-(0)=\chi_-\inv(g_-(\gamma(0)))=I$ we get
\begin{equation} \label{eqA}
\xi\circ\gamma=\chi_+.\xi,
\end{equation}
where `.' denotes the dressing action.

Moreover, we already know that
\begin{equation} \label{eqB}
(g_-\circ\gamma)(0)=\chi_-
\end{equation}
and
\begin{equation} \label{eqC}
\chi=\chi_-\chi_+\in\Lambda\LieSU(2)_\sigma.
\end{equation}
Now let us assume, that $\xi$ is of the form~\bref{xifconst}. 
Let $C=ce^{i\phi}$. Dressing with the matrix 
$\diag(\exp(-i\frac{\phi}2),
\exp(i\frac{\phi}2))\in\Lambda_r^+\LieSL(2,\C)_\sigma$ transforms $f$
into a positive real constant. On the other hand, 
by~\cite[Corollary~4.2]{DoHa:2}, dressing with a $\lambda$-independent
unitary matrix just amounts to a rigid rotation of the CMC immersion
in space. Thus w.l.o.g.\ we can assume that $C$ is a positive real constant.

Equation~\bref{eqA} gives
\begin{eqnarray}
\chi_+.f & = & (f\circ\gamma)\gamma'=\frac{C}{(\bquer z+\aquer)^2} \nonumber\\
& = & \frac{C\aquer^{-2}}{(1+\frac{\bquer}{C\aquer}(Cz))^2}=
\frac{C\aquer^{-2}}{(1+\frac{\bquer}{C\aquer}\int_0^z f\diff z)^2} \nonumber\\
& = & T_D(\sqrt{C}\aquer\inv)T_U(\frac{\bquer}{C\aquer})(f) \label{autodress}
\end{eqnarray}
where $f\equiv C$ and $T_U$, $T_D$ denote the basic dressing
transformations investigated in~\cite[Section~3]{DoHa:1}. I.e.\
$T_D(t)$ denotes dressing with $\diag(t,t\inv)$ and $T_U(t)$ denotes
dressing with $\tmatrix10{t\lambda}1$. 

Now Theorem~\ref{mainresult} implies that on surfaces with umbilics
the dressing action is free. 
Thus, \bref{autodress} determines the matrix $\chi_+$ uniquely:
\begin{equation}
\chi_+=\tmatrix{\frac{s}{\aquer}}00{\frac{\aquer}{s}}
\tmatrix10{\frac{\bquer\lambda}{C\aquer}}1=
\tmatrix{\frac{s}{\aquer}}0{\frac{\bquer\lambda}{Cs}}{
\frac{\aquer}{s}},
\end{equation}
where we have chosen $s=\sqrt{C}$ to be the positive root.
Now we have to find $\chi_-$ s.t.\ \bref{eqC} holds.
We make the ansatz
$$
\chi_-=\tmatrix1{q\lambda\inv}01.
$$
Then
$$
\chi=\chi_-\chi_+=\tmatrix{\frac{s}{\aquer}+\frac{\bquer q}{Cs}}{
\frac{q\aquer}{s}\lambda\inv}{\frac{\bquer}{Cs}\lambda}{\frac{\aquer}s}
$$
has to be unitary for all $\lambda\in S^1$. This is equivalent to
\begin{eqnarray} \label{eq3}
q = -\frac{b}{\aquer C}, \\
\frac{s}{\aquer}+\frac{\bquer q}{Cs} = \frac as.\label{eq4}
\end{eqnarray}
From this it follows, using $s^2=C$, that 
\begin{equation} \label{eq5prime}
C-\frac{|b|^2}{C^2}=|a|^2.
\end{equation}
The $\lambda\inv$-coefficient of Eq.~\bref{eqB} gives
$$
q=\int_0^{\gamma(0)}f\diff z=C\gamma(0)=C\frac b\aquer,
$$
which together with~\bref{eq3} gives $C=1$.
Finally,~\bref{eq5prime} implies for $C=1$, that $|a|^2+|b^2|=1$,
which together with $|a|^2-|b|^2=1$ forces $b=0$ and $|a|=1$, i.e.\
$\gamma$ is a rotation around the origin $z=0$ in $D_1$.
\QED

\newsection \label{nosyminautom}
The next best candidates for symmetric surfaces are expected to lie in
the dressing orbit of automorphic meromorphic potentials.
I.e. if $\gamma\in\Gamma$ is an automorphism in the symmetry group 
$\Aut_\Psi\FD$, $\FD=\C$ or $\FD=D_1$, then we look at a meromorphic
potential
$$
\xi_0=\lambda\inv\tmatrix0f{\frac{E}{f}}0\diff z,
$$
s.t.\ $E\circ\gamma={\gamma'}^{-2}E$, $f\circ\gamma={\gamma'}\inv f$.
Since $\xi_0$ is an automorphic one form w.r.t.\ $\gamma$, we
conclude for the integral $g_-^0$:
\begin{equation}
g_-^0(\gamma(z),\lambda)=\rho_-^0(\lambda)g_-^0(z,\lambda),
\end{equation}
where $\rho_-^0(\lambda)\in\Lambda^-_{r,\ast}\LieSL(2,\C)_\sigma$.

Let now $h_+\in\Lambda_r^+\LieSL(2,\C)_\sigma$ and define $g_-$ by the
dressing action of $h_+$ on $g_-^0$, 
\begin{equation} \label{gminfromg0}
h_+(\lambda)g_-^0(z,\lambda)=g_-(z,\lambda)p_+(z,\lambda)
\end{equation}
with $p_+:\FD\rightarrow\Lambda^+_r\LieSL(2,\C)_\sigma$. We denote by
$\xi$ the corresponding meromorphic potential $\xi=g_-\inv\diff g_-$.
We then get the following result:

\myprop{} {\em
Let $\gamma$ be an automorphism of $\FD$, $\FD=\C$ or $\FD=D_1$, and
let $\xi_0$ be a meromorphic potential with umbilics
on $\FD$ which satisfies $\gamma^\ast\xi_0=\xi_0$.
Let $h_+\in\Lambda_r^+\LieSL(2,\C)_\sigma$ and define $\xi$, $g_-^0$,
$\rho_-^0$, and $g_-$ as above. Then there is a unique
$\rho\in\Lambda_r\LieSL(2,\C)_\sigma$ and a unique map
$w_+:\FD\rightarrow\Lambda^+_r\LieSL(2,\C)_\sigma$, such that
\begin{equation} \label{gmintrafo}
g_-\circ\gamma=\rho g_-w_+.
\end{equation}
Moreover,
\begin{equation} \label{rhoeq}
\rho=h_+\rho_-^0h_+\inv.
\end{equation}
}

\Proof
First we note that
\begin{equation}
g_-\circ\gamma=h_+(g_-^0\circ\gamma)(p_+\inv\circ\gamma)
=h_+\rho_-^0h_+\inv g_-p_+(p_+\inv\circ\gamma).
\end{equation}
Thus,~\bref{gmintrafo} holds, with $w_+=p_+(p_+\inv\circ\gamma)$ and
$\rho=h_+\rho_-^0h_+\inv\in\Lambda_r\LieSL(2,\C)_\sigma$.

Now let $\rho\in\Lambda_r\LieSL(2,\C)_\sigma$ be an arbitrary matrix,
such that~\bref{gmintrafo} holds for some
$w_+:\FD\rightarrow\Lambda^+_r\LieSL(2,\C)_\sigma$.
Using the definition~\bref{gminfromg0} of $g_-$ we get from~\bref{gmintrafo}
\begin{equation} \label{isodresseq}
\rho\inv h_+\rho_-^0h_+\inv g_-=g_-w_+(p_+\circ\gamma)p_+\inv.
\end{equation}
Since $g_-(0,\lambda)=I$, we get
$\rho\inv h_+\rho_-^0h_+\inv\in\Lambda^+_r\LieSL(2,\C)_\sigma$.
Thus, $\rho\inv h_+\rho_-^0h_+\inv$ is in the isotropy group of
$\xi$ under dressing. 
Since we assumed that $\xi_0$ and therefore also $\xi$ has
umbilics, Theorem~\ref{mainresult} gives
\begin{equation}
\rho=h_+\rho_-^0h_+\inv,
\end{equation}
from which the uniqueness of $w_+$ also follows.
\QED

\mytheorem{} {\em
Under the same assumptions as in the proposition above 
the following two statements are equivalent:
\begin{enumerate}
\item The automorphism $\gamma$ is in the symmetry group of the
CMC immersion $\Psi$ corresponding to $\xi$, i.e.\ $\gamma\in\Aut_\Psi\FD$,
\item The integral $g_-$ of $\xi$ is invariant under $\gamma$,
i.e.\ $g_-\circ\gamma=g_-$.
\end{enumerate}
}

\Proof
$2.\Rightarrow 1.$ follows immediately from \cite[Corollary~4.2]{DoHa:2}.
For the converse statement assume that
$\gamma\in\Aut_\Psi\FD$. Then $g_-$ transforms by
\cite[Theorem~4.2]{DoHa:2} and~\cite[Lemma~2.2]{DoHa:3} like~\bref{gmintrafo},
where $\rho\in\Lambda_r\LieSU(2)_\sigma\subset\Lambda_r\LieSL(2,\C)_\sigma$ and
$w_+:\FD\rightarrow\Lambda^+_r\LieSL(2,\C)_\sigma$. 
By Proposition~\ref{nosyminautom} we have $\rho=h_+\rho_-^0h_+\inv$.
Hence, on the circle $S_r$ the eigenvalues of $\rho$ and $\rho_-^0$ 
coincide. The matrix $\rho$ is by \cite[Section~3]{DoHa:2} also the monodromy
matrix of the extended frame $F$ of $\Psi$,
\begin{equation}
F(\gamma(z),\lambda)=\rho(\lambda)F(z,\lambda)k(z)
\end{equation}
for some $k:\FD\rightarrow\LieU(1)$. Therefore,
by \cite[Lemma~2.2]{DoHa:3}, the unitary matrix $\rho$ can be extended
holomorphically to $\C^\ast$. And since $\rho_-^0$
can be extended holomorphically to the exterior of the circle $S_r$,
we get that the eigenvalues of $\rho$ and $\rho_-^0$ are holomorphic
functions on $\CPE$. Thus,
they are constant and equal to the eigenvalue of $\rho_-^0$ at
$\lambda=\infty$. By $\rho_-^0(\lambda\rightarrow\infty)=I$ we get
$\rho=I$, which together with \bref{gmintrafo} gives the desired result.
\QED

\myremark{}
It should be noted here, that due to Theorem~\ref{nosyminautom}, the
situation for CMC immersions with umbilics is very different from
the situation without umbilics: As was shown
in~\cite{DoWu:1}, {\em all} CMC immersions of finite type, among which
are all CMC tori, can be obtained by dressing the translationally
invariant meromorphic potential $\xi=\lambda\inv\tmatrix0110\diff z$ of the
standard cylinder. However, as follows from~\cite[Eq.~(3.4.5)]{DoHa:3}, {\em
none} of these surfaces (not even the standard cylinder itself)
has a constant monodromy matrix $\rho=\chi\equiv I$. I.e.\
Theorem~\ref{nosyminautom} can obviously not be extended to potentials
without umbilics. 
Therefore, there is no immediate way to generalize the application of
the dressing group to finite type surfaces, as it was done
in~\cite{DoHa:3}, to surfaces with umbilics. In fact the triviality of
the dressing isotropy groups of CMC immersions with umbilics seems to
indicate, that for surfaces with umbilics, in particular compact
surfaces of genus $g\geq2$, there is no algebraic geometric method
corresponding to Krichever's method.


\begin{thebibliography}{10}

\bibitem{AdvMo:1}
{\sc M.~Adler and P.~{van Moerbeke}}, {\em Completely integrable systems,
  euclidean {L}ie algebras and curves}, Adv. in Math., 38 (1980), pp.~318--379.

\bibitem{Bo:1}
{\sc A.~Bobenko}, {\em All constant mean curvature tori in ${R}^3$, ${S}^3$,
  ${H}^3$ in terms of theta-functions}, Math. Ann., 290 (1991), pp.~209--245.

\bibitem{BuPe:1}
{\sc F.~Burstall and F.~Pedit}, {\em Dressing {O}rbits of {H}armonic {M}aps},
  Duke Math. J., 80 (1995), pp.~353--382.

\bibitem{DoHa:3}
{\sc J.~Dorfmeister and G.~Haak}, {\em On constant mean curvature surfaces with
  periodic metric}, Pac.J.Math.,  (1996). dg-ga/9609006
\newblock to appear.

\bibitem{DoHa:2}
\leavevmode\vrule height 2pt depth -1.6pt width 23pt, {\em On symmetries of
  constant mean curvature surfaces}.
\newblock Sfb 288 preprint 197, 1996. dg-ga/9603007

\bibitem{DoHa:1}
\leavevmode\vrule height 2pt depth -1.6pt width 23pt, {\em Meromorphic
  potentials and smooth surfaces of constant mean curvature}, Math. Z., 224
  (1997), pp.~603--640. dg-ga/9412007

\bibitem{DoPeWu:1}
{\sc J.~Dorfmeister, F.~Pedit, and H.~Wu}, {\em Weierstra\ss\ {T}ype
  {R}epresentations of {H}armonic {M}aps into {S}ymmetric spaces}, Comm.
  Analysis and Geom.,  (1996).
\newblock to appear.

\bibitem{DoWu:1}
{\sc J.~Dorfmeister and H.~Wu}, {\em Constant mean curvature surfaces and loop
  groups}, J. reine angew. Math., 440 (1993), pp.~43--76.

\bibitem{Mc:1}
{\sc I.~McIntosh}, {\em Global solutions of the elliptic 2{D} periodic {T}oda
  lattice}, Nonlinearity, 7 (1994), pp.~85--108.

\bibitem{PiSt:1}
{\sc U.~Pinkall and I.~Sterling}, {\em On the classification of constant mean
  curvature tori}, Annals of Math., 130 (1989), pp.~407--451.

\bibitem{ReSe:1}
{\sc A.~Reyman and M.A.Semenov-Tian-Shansky}, {\em Reduction of {H}amiltonian
  systems, affine {L}ie algebras and {L}ax equations}, Invent. Math., 54
  (1979), pp.~81--100.

\bibitem{SeWi:1}
{\sc G.~Segal and G.~Wilson}, {\em Loop groups and equations of {K}d{V} type},
  Publ. Math. I.H.E.S., 61 (1985), pp.~5--65.

\bibitem{Wu:3}
{\sc H.~Wu}, {\em Denseness of plain constant mean curvature surfaces}.
\newblock preprint, 1996.

\bibitem{Wu:4}
\leavevmode\vrule height 2pt depth -1.6pt width 23pt, {\em On the dressing
  action of loop groups on constant mean curvature surfaces}.
\newblock preprint, 1996.

\bibitem{Wu:5}
\leavevmode\vrule height 2pt depth -1.6pt width 23pt, {\em A new
  characterization of normalized potentials in dimension two}.
\newblock preprint, 1997.

\bibitem{ZaSh:1}
{\sc V.~Zakharov and A.~Shabat}, {\em Integration of the nonlinear
  {S}chr\"odinger equations of mathematical physics by method of inverse
  scattering problem,~ii ({E}ngl. transl.)}, Funct. Anal. Appl., 13 (1979),
  pp.~166--174.

\end{thebibliography}
\end{document}